# The dunce hat in a minimal non-extendably collapsible 3-ball

January 28, 2013

### Authors

Bruno Benedetti and Frank H. Lutz

### Description

We obtain a geometric realization of a minimal 8-vertex triangulation of the dunce hat in Euclidean 3-space. We show there is a simplicial 3-ball with 8 vertices that is collapsible, but also collapses onto the dunce hat, which is not collapsible. This 3-ball is as small as possible, because all triangulated 3-balls with fewer vertices are extendably collapsible. As we will see, the Alexander dual of the dunce hat is collapsible.

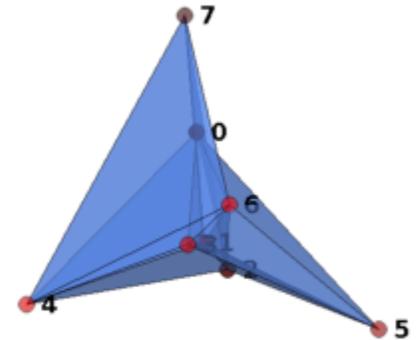

Figure 1: The dunce hat in 3-space.

The *dunce hat* [29] is the 2-dimensional cell complex obtained from one triangle by identifying all three boundary edges, provided not all edges are oriented coherently. The dunce hat can be triangulated minimally with 8 vertices in several ways. We show that one of these triangulations D sits in the boundary of a simplicial 4-polytope with precisely 8 vertices. In fact, the polytope is No. 32 in the list by Grünbaum and Sreedharan [16]. Via a Schlegel diagram, this yields a geometric realization of this 8-vertex dunce hat D in Euclidean 3-space. The boundary of the 4-polytope No. 32 is a simplicial 3-sphere with 19 tetrahedra. Seven of them can be removed to yield a 3-ball B with 8 vertices that is collapsible but not extendably collapsible, as D is still contained in the 2-skeleton of B. This B is even vertex-decomposable, shellable and non-evasive. As a consequence, we show that the Alexander dual of the dunce hat is collapsible, thereby clarifying an argument by Kahn, Saks and Sturtevant [21; p. 301].

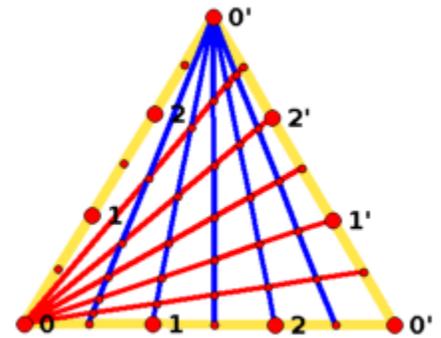

Figure 2: The dunce hat; Step I.

Let us recall some basic definitions. A face f of a simplicial complex K is *free* if it is properly contained in only one other face F of K. An *elementary collapse* is the simultaneous removal from a simplicial complex of a free face f together with the unique face F containing f. Topologically, every elementary collapse is a strong deformation retract: It corresponds to pushing f inside F. A simplicial complex K is *collapsible* if some sequence of elementary collapses reduces it to a single point. A complex K is *extendably collapsible* if every initial sequence (possibly empty) of elementary collapses can be extended to a sequence that reduces K to a point. More generally, one says that K *collapses onto* H to mean that a complex K can be reduced onto a subcomplex H via some sequence of elementary collapses.

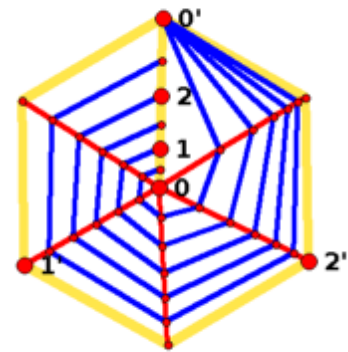

Figure 3: The dunce hat; Step II.

Since elementary collapses preserve the homotopy type, all collapsible complexes are contractible. For 1-dimensional complexes the two notions are equivalent. The contractible 1-complexes are trees and therefore are collapsible, by removing one leaf at the time. In contrast, the dunce hat is a contractible 2-dimensional space none of whose triangulations is collapsible (as there are no free edges to begin with). It is observed in [3] that all $Z_2$-acyclic (in particular, all contractible) simplicial complexes with up to 7 vertices are collapsible. Thus, at least 8 vertices are needed to triangulate the dunce hat. Triangulations with 8 vertices are easy to construct; see for example Figures 6 and 7. There are also other examples of non-collapsible contractible complexes, like Bing's house with two rooms [7] or the abalone [20]. Indeed, a simplicial complex K is contractible if and only if there is an integer d such that (the barycentric subdivision of) $K \times I^d$ is collapsible [1], [11].

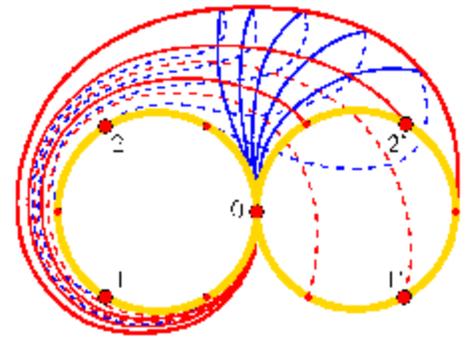

Figure 4: The dunce hat; Step III.

In the following, $S_{32}$ will denote the boundary of the simplicial 4-polytope No. 32 in the list by Grünbaum-Sreedharan. The 3-sphere $S_{32}$ has 8 vertices and 19 facets.

**Theorem 1:** The 8-vertex sphere $S_{32}$ contains an 8-vertex triangulation of the dunce hat as a subcomplex. In particular, a Schlegel diagram of $S_{32}$ realizes the dunce hat geometrically in Euclidean 3-space.

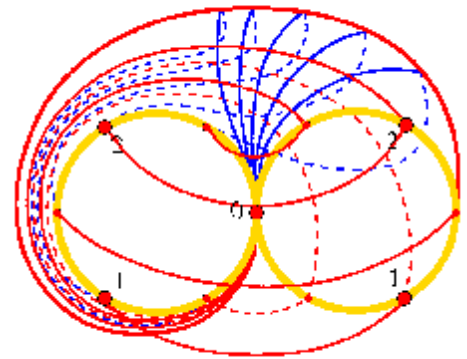

Figure 5: The dunce hat; Step IV.

*Proof:* In a triangle that has not all three boundary edges oriented coherently, there are two choices for a pair of oppositely oriented edges; glueing together such two edges yields a cone. In Figure 2, the red and the blue fans indicate the two possible cones, which we call the *red cone* and the *blue cone*, respectively. Once we have formed, say, the red cone (cf. Figure 3), we still have to glue together the yellow seam-ray 0-1-2-0' and the yellow boundary circle 0'-1'-2'-0'. This glueing is depicted in Figures 4-5. Here, we first identify the apex 0 with the vertex 0' in the boundary circle (thus closing the yellow seam-ray 0-1-2-0' to a loop 0-1-2-0, cf. Figure 4) and then we glue together the two circles (as displayed in Figure 5). We apply this construction to the 8-vertex triangulation D of the dunce hat from Figure 6 and Figure 7. If we add to the red cone of Figure 6 the triangle 0'1'2', this encloses a solid cone. In the same manner, we can add the triangle 0"1"2" to the blue cone of Figure 7 to encompass again a solid cone. We then glue together the two solid cones along their common boundary. The manifold we get is a 3-sphere. Next, we triangulate the interiors of the two solids. We cut the red solid into three pieces by introducing the spanning triangles 134 and 067. The purpose of introducing these two triangles is to separate identified vertices so that we will not get multiple edges in the interior. The tetrahedron 0134 fills the first piece (next to the apex), while the two tetrahedra

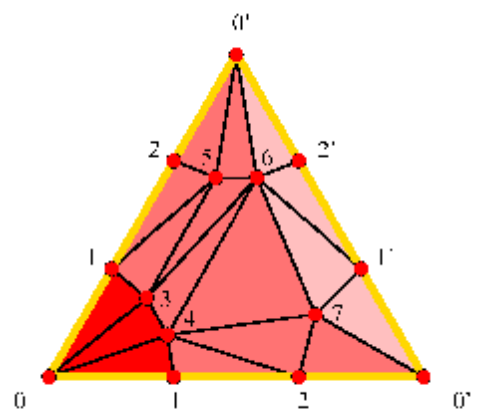

Figure 6: The dunce hat, red cone.

0126 and 0167 fill the bottom piece. For the central piece we use seven tetrahedra (0257, 0567, 1245, 1345, 2457, 3456, 4567). Likewise, the blue cone is split into three pieces by introducing the two spanning triangles 256 and 047. The three pieces of the blue solid are then filled with the tetrahedron 0256 at the apex, with the six tetrahedra 0137, 0347, 1256, 1356, 1367, 3467 in the central piece, and with 0124 and 0247 in the third piece. Altogether, these 19 tetrahedra yield a triangulated 3-sphere with 8 vertices. This sphere is combinatorially isomorphic to $S_{32}$. Coordinates for $S_{32}$ are obtained using the `rand_sphere` client of `polymake` [15]. Via a Schlegel diagram, this yields a geometric realization of D in Euclidean 3-space; see Figure 1. Q.e.d.

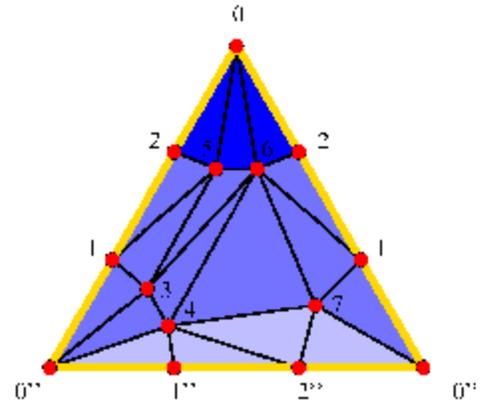

Figure 7: The dunce hat, blue cone.

There are exactly 39 triangulated 3-manifolds with 8 vertices. All of them are triangulated spheres [2], and 37 of them are polytopal. The triangulation D of the dunce hat is also contained in Grünbaum-Sreedharan's sphere No. 37 (polytopal, with 8 vertices and 20 facets) and in Grünbaum-Sreedharan's sphere M (non-polytopal, with 8 vertices and 20 facets) [16], but it is not contained in any other of the 39 3-spheres with 8 vertices, as we checked by computer.

We next construct a ball that is collapsible, but not extendably collapsible. The ball is even vertex-decomposable, shellable and non-evasive, where the first property implies the other two. Shellable contractible complexes are always collapsible: The facets can be collapsed according to the reverse of a shelling order. On the contrary, extendably shellable complexes [10] need not be extendably collapsible, as we will see in Theorem 2 below.

**Theorem 2:** Seven tetrahedra can be removed from $S_{32}$ to yield a shellable 3-ball B, which is collapsible, extendably shellable, but not extendably collapsible. Any ball with fewer vertices than B is extendably collapsible.

*Proof:* Removing the seven tetrahedra 0257, 0567, 1245, 1345, 2457, 3456, 4567 (the central red piece) from the sphere No. 32, we obtain a 3-ball B with 8 vertices and 12 tetrahedra that still contains the dunce hat D as a subcomplex. We used the classification of all 3-balls with up to 9 vertices from [23] to show that no triangulated 3-ball on 8 vertices with less than 12 tetrahedra contains D. Furthermore, B is the only triangulated 3-ball with 8 vertices and 12 tetrahedra that contains D. Note that all 3-balls with up to 8 vertices are (extendably) shellable and therefore collapsible. However, for a suitable choice of collapsing steps, the ball B collapses onto the dunce hat D in its 2-skeleton. Thus, B is not extendably collapsible. By [3] all contractible 2-complexes with up to 7 vertices are collapsible and it follows that all 3-balls with up to 7 vertices are extendably collapsible. Q.e.d.

Any triangulation of the dunce hat is non-shellable, not constructible [17], but Cohen-Macaulay over any field (as can be seen by using Reisner's criterion [27]). It is not known whether all constructible contractible complexes are collapsible. However, all constructible 3-balls are collapsible [6]. (In contrast, collapsible balls can be non-constructible. Lickorish showed that a collapsible ball might contain a 3-edge knot in its 1-skeleton [22], while a 3-edge closed path inside a shellable or a constructible 3-ball is always unknotted [18].)

**Theorem 3:** The ball B is vertex-decomposable and non-evasive. In particular, the Alexander dual of the dunce hat is evasive but collapsible.

*Proof:* A valid sequence of vertex deletions for B is 5-2-4-0-1. Vertex-decomposability implies non-evasiveness for balls [5]. Non-evasive complexes are collapsible [21]; hence the 8-vertex dunce hat D is evasive. It is proven in [21] that a simplicial complex is evasive if and only if its Alexander dual is evasive. Therefore the Alexander dual of D is evasive, while the Alexander dual of B is non-evasive (and in particular collapsible). As explained in [21], since B collapses onto D, the Alexander dual of D collapses onto the Alexander dual of B, which is collapsible. Q.e.d.

In the 3-sphere described in Theorem 1, the complement of a regular neighborhood of the dunce hat D is topologically a 3-ball. However, there are embeddings of the dunce hat in $S^4$ so that the complements of the respective regular neighborhoods are not even simply connected [25]. By a result of Freed [13], the dunce hat has infinitely many distinct embeddings into $S^4$, distinguished by the fundamental groups of the respective complements. Two dunce hats can link in $S^4$ [29].

All collapsible 1- and 2-dimensional complexes are extendably collapsible. In fact, a graph is collapsible if and only if it is a tree. If instead K is 2-dimensional and contains a pure 2-complex C without free faces, it is easy to see that K is not collapsible either. In contrast, in dimension 3, if we start with the 3-ball B of Theorem 2, some collapsing sequences get stuck and some do not. Moreover, some simplicial subdivision of the product of the dunce hat with an interval collapses onto a point as well as onto (a subdivision of) the dunce hat [29]. In view of this result, Zeeman conjectured that for *every* contractible 2-complex K the product of K with an interval I has a collapsible subdivision [29]. Zeeman's conjecture (which is still open) implies the 3-dimensional Poincaré conjecture (as proved by Perelman [26]) and the Andrews-Curtis conjecture (still open, cf. [19]). For recent developments see [1].

Cohen proved (by assuming that the 4-dimensional Poincaré conjecture is true, as later proved by Freedman [14]) that for every contractible 2-complex K the product of K with $I^5$ has a collapsible subdivision [8]. More generally, for any contractible complex K the product of K with a sufficiently high dimensional cube has a collapsible subdivision [8]. However, for every d>2, there is a contractible d-complex such that the product with a single interval has no collapsible subdivision [9], [19].

If a contractible 2-complex K embeds into a 3-manifold, then K x $I^3$ has a collapsible subdivision [8], [29]. However, some contractible 2-complexes do not embed into any 3-manifold, for example, cones over non-planar graphs [29]. For d = 3 and d ≥ 5, a PL d-manifold that collapses onto a contractible 2-complex is a PL d-ball [29], but Zeeman [29] showed that the Mazur 4-manifold [24], although different from a 4-ball, collapses onto the dunce hat. No triangulation of the Mazur manifold is collapsible (because all collapsible PL triangulations are PL balls [28]).

A sequence of elementary collapses can always be rearranged so that higher-dimensional faces are collapsed first. The initial segment of a collapsing sequence is particularly nice to describe in case the starting complex is a manifold M minus a facet F. In fact, in every collapse of M-F onto a lower-dimensional complex, the d-faces of M are removed along a (directed) spanning tree T of the dual graph of M, rooted at F. As a result, M-F collapses onto the (d-1)-complex $K^T$ of all the (d-1)-faces that are not perforated by T. In particular, M-F is collapsible if and only if $K^T$ is collapsible for some dual spanning tree T [6]. This $K^T$ does not depend on F: Thus, either M-F collapses for all facets F, or M-F does not collapse for any F. However, the collapsibility of $K^T$ does depend on T. An example is given by the 3-sphere $S_{32}$ presented in Theorem 1. In fact, for any tetrahedron F, the 3-ball $S_{32}$-F may collapse onto a point or also onto the dunce hat D, depending on which tree we select. This has recently been connected in [4] and [6] with the notion of local constructibility, introduced in the quantum gravity paper by Durhuus and Jónsson [12].

## Acknowledgment

The first author was supported by the DFG via the Berlin Mathematical School, the second author was supported by the DFG Research Group "Polyhedral Surfaces", Berlin.

| | |
|---|---|
| **Keywords** | dunce hat; 3-ball; vertex-minimal triangulation; non-extendably collapsible; collapsibility and contractibility of simplicial complexes; shellability |
| **MSC-2000 Classification** | 57Q05 (57Q15, 57Q35) |

**Authors' Addresses**

Bruno Benedetti

    Royal Institute of Technology (KTH)
    Department of Mathematics
    S-100 44 Stockholm
    Sweden
    brunoben@math.kth.se
    http://www.math.kth.se/~brunoben/

Frank H. Lutz

    Technische Universität Berlin
    Institut für Mathematik, MA 3-2
    10623 Berlin
    Germany
    lutz@math.tu-berlin.de
    http://page.math.tu-berlin.de/~lutz/